\documentclass[11pt]{amsart}
\begin{document}

\title{On the total mean curvature of non-rigid surfaces}
\author{Victor Alexandrov}
\thanks{The author was supported in part by the Russian State Program for Leading Scientific Schools, Grant~NSh--8526.2008.1.}
\address{Sobolev Institute of Mathematics, Novosibirsk, 630090, Russia and Novosibirsk State University, Novosibirsk, 630090, Russia}
\email{alex@math.nsc.ru}
\date{November 29, 2008}
\begin{abstract}
Using Green's theorem we reduce the variation of the total mean curvature of a smooth
surface in the Euclidean 3-space to a line integral of a special vector field
and obtain the following well-known theorem as an immediate consequence:
the total mean curvature of a closed smooth surface in 
the Euclidean 3-space is stationary under an infinitesimal flex.
\end{abstract}
\maketitle

A smooth surface $S\subset\mathbb R^3$ is said to be \textit{flexible}
if there is a smooth mapping $\boldsymbol{\varphi}: S\times (-1,1)\to\mathbb R^3$
such that

(1) for every smooth curve $\gamma\subset S$ and every $t\in (-1,1)$
the length of the curve $\{\boldsymbol{\varphi}(\boldsymbol{x},t)\vert \boldsymbol{x}\in\gamma\}$ is
equal to the length of $\gamma$;

(2) for every $t\neq r$ there are two points $\boldsymbol{x},\boldsymbol{y}\in S$ such that
the Euclidean distance between the points $\boldsymbol{\varphi}(\boldsymbol{x},t)$ 
and $\boldsymbol{\varphi}(\boldsymbol{y},t)$ is not equal to the Euclidean
distance between the points $\boldsymbol{\varphi}(\boldsymbol{x},r)$ and 
$\boldsymbol{\varphi}(\boldsymbol{y},r)$.

In other words, a smooth surface $S\subset\mathbb R^3$ is said to be flexible
if there exists a family $\{S_t\}_{t\in (-1,1)}$ of smooth surfaces $S_t\subset\mathbb R^3$,
such that (a)~$S_0=S$; (b)~$S_t$ is isometric to $S_0$ in the intrinsic metrics (see, e.\,g., [\textbf{1}] 
for detail) for every $t$; and (c)~$S_t$ and $S_r$ are not congruent if $t\neq r$.

One can easily check that a plane disk is flexible, but a very long-standing problem
reads that no compact boundary-free smooth surface in $\mathbb R^3$ is flexible [\textbf{8}; Problem 50].
The reader interested in a similar problem for polyhedral surfaces is referred to [\textbf{5}].

If $S$ is oriented then the \textit{total mean curvature} of $S_t$ is given by the classical formula
$$
H(S_t)=\int_{S_t}\frac{1}{2}\bigl(\kappa_1(\boldsymbol{x})+\kappa_2(\boldsymbol{x})\bigr)\,d(S_t)\eqno(1)
$$
where $\kappa_1(\boldsymbol{x})$ and $\kappa_2(\boldsymbol{x})$ are the principal curvatures of $S_t$
at the surface point $\boldsymbol{x}$.

A smooth surface $S\subset\mathbb R^3$ is said to be \textit{non-rigid}
if there is a smooth vector field $\boldsymbol{v}: S\to\mathbb R^3$
such that

(i) the smooth mapping $\boldsymbol{\psi}: S\times (-1,1)\to\mathbb R^3$
defined by the formula $\boldsymbol{x} \mapsto \boldsymbol{x}+t\boldsymbol{v}(\boldsymbol{x})$ is such that,
for every smooth curve $\gamma\subset S$, the length of the curve 
$\{\boldsymbol{\psi}(\boldsymbol{x},t)\vert \boldsymbol{x}\in\gamma\}$ is
stationary at $t=0$;

(ii) no family of rigid motions of $S$ generates $\boldsymbol{v}$.

The above field $\boldsymbol{v}$ is called a (non-trivial) \textit{infinitesimal flex} of $S$.
Non-rigid compact surfaces in~ $\mathbb R^3$ do exist and were studied by many authors
(references may be found, e.\,g., in [\textbf{4}]). In particular,
it is known that if~ $S$ is parameterized by $\boldsymbol{x}=\boldsymbol{x}(u,v)$ 
and $\boldsymbol{v}=\boldsymbol{v}(u,v)$ is its infinitesimal flex then 
$$\boldsymbol{x}_u\cdot\boldsymbol{v}_u=\boldsymbol{0},\quad
\boldsymbol{x}_u\cdot\boldsymbol{v}_v+\boldsymbol{x}_v\cdot\boldsymbol{v}_u=\boldsymbol{0},\quad\text{and}\quad
\boldsymbol{x}_v\cdot\boldsymbol{v}_v=\boldsymbol{0},\eqno(2)
$$
where $\cdot$ stands for the scalar product in $\mathbb R^3$.

Let $S$ be a compact oriented smooth surface in $\mathbb R^3$.
Note that, for all~ $t$ close enough to zero, the surface 
$\boldsymbol{\psi} (S,t)=\{\boldsymbol{\psi}(\boldsymbol{x},t)\vert \boldsymbol{x}\in S\}$
is smooth and oriented. Denote by $\boldsymbol{n}(\boldsymbol{x}, t)$ its unit normal vector
to the surface $\boldsymbol{\psi} (S,t)$ at the point $\boldsymbol{\psi}(\boldsymbol{x},t)$
and denote by $\boldsymbol{n'}(\boldsymbol{x},t)$ the velocity vector of the vector-function 
$t\mapsto\boldsymbol{n}(\boldsymbol{x},t)$, i.\,e., put by definition
$$
\boldsymbol{n'}(\boldsymbol{x},t)=\frac{d}{dt}\boldsymbol{n}(\boldsymbol{x},t).
$$
Define the vector field $\boldsymbol{m}$ on $S$ by the formula
$\boldsymbol{m}(\boldsymbol{x})=\boldsymbol{n'}(\boldsymbol{x},0)\times\boldsymbol{n}(\boldsymbol{x},0)$,
where $\times$ stands for the cross product in $\mathbb R^3$.
(Note that $\boldsymbol{m}$ is a tangential vector field on $S$, though we will
not use this fact below.)
At last, put by definition
$$
H'(S)=\frac{d}{dt}\biggl|_{t=0}H\bigl(\boldsymbol{\psi} (S,t)\bigr).
$$
For obvious reasons, we call $H'(S)$ the \textit{variation} 
of the total mean curvature of $S$.

The main result of this note reads as follows:

\textbf{Theorem.} \textit{For every compact oriented smooth surface $S$ in $\mathbb R^3$
and any its infinitesimal flex $\boldsymbol{v}$, the variation of the total mean curvature of ~$S$
equals the line integral of the vector field  $\boldsymbol{m}$ over the boundary $\partial S$
of $S$, i.\,e.,}
$$
H'(S)=\frac{1}{2}\int_{\partial S}\boldsymbol{m}(\boldsymbol{x})\cdot d\boldsymbol{x}.
$$
(Of course, the curve $\partial S$ of the line integral is supposed to have positive orientation.) 

\textit{Proof}. It suffice to prove the theorem `locally', i.\,e.,
for a surface $S$ covered with a single chart.
In particular, we may assume that $S$ is parameterized by 
$\boldsymbol{x}=\boldsymbol{x}(u,v)=\bigl(u,v,f(u,v)\bigr)$,
$(u,v)\in D\subset R^2$. Let
$\boldsymbol{v}(\boldsymbol{x})=\boldsymbol{v}(u,v)=\bigl(\xi(u,v),\eta(u,v),\zeta(u,v)\bigr)$.
Then equations (2) take the form
$$
\begin{cases}
\xi_u=-f_u\zeta_u,\\ 
\xi_v+\eta_u=-f_v\zeta_u-f_u\zeta_v,\\ 
\eta_v=-f_v\zeta_v. 
\end{cases}\eqno(3)
$$
Differentiating (3) with respect to $u$ and $v$ and calculating
linear combinations yields
$$
\begin{cases}
\xi_{uu}=-f_{uu}\zeta_u-f_u\zeta_{uu},\\
\xi_{uv}=-f_{uv}\zeta_u-f_u\zeta_{uv},\\
\xi_{vv}=-f_{vv}\zeta_u-f_u\zeta_{vv},\\
\eta_{uu}=-f_{uu}\zeta_v-f_v\zeta_{uu},\\
\eta_{uv}=-f_{uv}\zeta_v-f_v\zeta_{uv},\\
\eta_{vv}=-f_{vv}\zeta_v-f_v\zeta_{vv}.
\end{cases}\eqno(4)
$$

Suppose $S$ is oriented by the following field of the unit normal vectors
$(1+f_u^2+f_v^2)^{-1/2}\bigl(-f_u,-f_v,1\bigr)$.
Using equations (3) and the standard machinery of differential geometry [\textbf{3}], we
get
$$
2H'(S)=\iint_D \bigl[
(1+f_v^2)\zeta_{uu}-2f_uf_v\zeta_{uv}+(1+f_u^2)\zeta_{vv}
\bigr]\,dudv.\eqno(5)
$$

On the other hand, direct calculations show that
\begin{multline*}
\boldsymbol{m}(u,v)=\frac{1}{1+f_u^2+f_v^2}
\bigl(f_u\xi_v+f_v\eta_v-\zeta_v,
-f_u\xi_u-f_v\eta_u+\zeta_u,\\
-(f_u^2+f_v^2)\eta_u+f_v\zeta_u-f_u(1+f_u^2+f_v^2)\zeta_v
\bigr)
\end{multline*}
and
$$
\int_{\partial S}\boldsymbol{m}(\boldsymbol{x})\cdot d\boldsymbol{x}
=\int_{\partial D}
\bigl((1+f_u^2)\zeta_v-f_u\eta_u\bigr) du +\bigl(\zeta_u-f_v\eta_u-f_uf_v\zeta_v\bigr) dv.\eqno(6)
$$
Applying Green's theorem 
$$\int_{\partial D}P\,du+Q\,dv=\iint_D\biggl(\frac{\partial Q}{\partial u}-\frac{\partial P}{\partial v}\biggr)dudv$$
to the right-hand size integral in (6) and using formulas (3) and (4), we
transform (6) to the right-hand size integral in (5).\hfill q.e.d.

\textbf{Corollary 1.} \textit{For every compact oriented boundary-free smooth surface $S$ in $\mathbb R^3$
and any its infinitesimal flex, the variation of the total mean curvature of $S$ equals zero.}

\textbf{Corollary 2.} \textit{Every flexible compact oriented boundary-free smooth surface in $\mathbb R^3$
preserves its total mean curvature during the flex.}

The both corollaries immediately follow from the above theorem. 
In fact, they are known in a much more general situation, 
namely, for piecewise smooth hypersurfaces in multidimensional Euclidean spaces
(see [\textbf{2}], [\textbf{6}], and [\textbf{7}]).
But reduction to a line integral is new and, probably, may help to understand what 
other quantities remain constant during the flex. 
Integrals of symmetric functions of the principal curvatures?
Volume? The reader interested in similar results for polyhedra
is referred to [\textbf{5}] and literature mentioned therein.

\end{document}